\documentclass[12pt]{amsart}
\usepackage{amsmath}
\usepackage{amssymb}
\usepackage{amsfonts}
\usepackage{amsthm}
\usepackage{verbatim}
\usepackage{amscd}
\usepackage{cite}
\usepackage{leftidx}
\usepackage{enumerate}
\usepackage{txfonts}
\usepackage{manfnt}
\usepackage{amscd}
\usepackage{mathpazo}
\usepackage[mathscr]{eucal}
\usepackage{hyperref}
\DeclareMathOperator{\mgbar}{\overline\M_g}

\DeclareMathOperator{\coker}{coker}

\def\refp #1.{(\ref{#1})}

\newcommand{\M}{\mathcal{M}}

\newcommand{\lbl}[1]{\label{#1}}

\newcommand{\Cal}[1]{\mathcal #1}

\def\sbr #1.{^{[#1]}}
\def\sfl #1.{^{\lfloor #1\rfloor}}

\newcommand{\myit}[1]{\emph{\ {#1}\ }}

\newcommand{\sectionend}
{\[\circ\gtrdot\ggg\mathbb{\times}\lll\lessdot\circ\]}

\newcommand\red{{\mathrm red}}

\DeclareMathOperator{\Mod}{Mod}

\def\?{{\bf{??}}}

\def\M{\Cal M}

\def\L{\mathcal L}

\def\O{\mathcal O}

\def\g{\mathfrak g}

\def\1/2{\frac{1}{2}}

\def\2{{[2]}}

\def\<{\langle}
\def\>{\rangle}

\def\2{{[2]}}

\def\scl #1.{^{\lceil#1\rceil}}
\def\spr #1.{^{(#1)}}
\def\sbc #1.{^{\{#1\}}}

\def\subpr#1.{_{(#1)}}

\def\beq{\begin{equation*}}
\def\eeq{\end{equation*}}

\newcommand{\lf}[1]{\leftidx{_{\ \mathrm L}}{#1}{}}
\newcommand{\rt}[1]{\leftidx{_{\ \mathrm R}}{#1}{}}

\def\g3{{\Gamma\spr 3.}}

\def\ggg{{\Gamma\spr 3.}}

\newcommand{\eqspl}[2]{
\begin{equation}\label{#1}
\begin{split}
#2\end{split}\end{equation}}
\newcommand{\eqsp}[1]{\begin{equation*}
\begin{split}#1\end{split}\end{equation*}}

\newcommand{\exseq}[3]{
0\to #1\to #2\to #3\to 0
}
\newcommand{\beginalphaenum}{
\begin{enumerate}\renewcommand{\labelenumi}{ }
\item \begin{enumerate}
}

\def\eex{\end{rm}\end{example}}
\newcommand\newsection[1]{\section{#1}\setcounter{equation}{0}
}

\newtheorem{thm}{Theorem}[section]

\newtheorem*{thm*}{Theorem}
\newtheorem{cor}[thm]{Corollary}
\newtheorem*{cor*}{Corollary}

\newtheorem*{claim*}{Claim}
\newtheorem{prop}[thm]{Proposition}

\newtheorem{defn}[thm]{Definition}%number together with conj
\theoremstyle{remark}

\newtheorem{rem}[thm]{Remark}

\newtheorem{example}[thm]{Example}
\newtheorem*{example*}{Example}

\begin{document}
%\reversemarginpar
%\Large
%{\fontfamily{duerer} \selectfont cdr12}
\title {Echelon modifications of vector bundles}

\author
{Ziv Ran}
%\thanks{\raggedright{
%Partially supported by NSA Grant MDA904-02-1-0094} }
\date {\today}% \enddate
%\affil {University of California, Riverside}\endaffil
\address{\tiny  {\newline Ziv Ran \newline University of California
Mathematics Department\newline Surge Facility, Big Springs Rd.
\newline Riverside CA 92521
\newline ziv.ran @ ucr.edu}}
%\email { ziv.ran @ ucr.edu}

%\urladdr http://www.math.ucr.edu/~ziv/papers/
%semicurv.pdf\endurladddr
%\rightheadtext { Cycle Map and Intersection Calculus}
%\leftheadtext{Z. Ran}
\begin{abstract}
%We compute, in terms of tautological classes, the fundamental class of
We study a filtered generalization of the operation of elementary modification of vector bundles. The generalization is motivated by applications to the degeneration theory of
 linear systems.

\end{abstract}

\maketitle
%\newsection{Echelon modifications}\lbl{echelon-mod}
The notion of  \emph{elementary modification}
of a vector
bundle along a divisor on a scheme is well known and is a standard method for constructing vector bundles.
The purpose of this note is to define a generalization of
this  to a certain \emph{filtered} setting,
 on both the bundle and divisor sides.
 The generalization, whose main properties are given in
  Theorem \ref{thm} below, is motivated by situations which occur
 in the study of linear systems on a family of curves
 with reducible fibres.
 The construction made here  will be
applied in \cite{grd} to yield boundary modifications of the Hodge bundle on the
moduli space of curves. These modifications will play a key role in
our work on the closure in $\mgbar$ of $g^r_d$ loci in $\M_g$. See also \cite{canonodal}.\par
We will work in a general setting of vector bundles on a scheme. However,
we wish to draw attention to a couple of relatively subtle points arising in the construction, which may seem
 surprising from such a general vantage point,
 and which are both motivated by the
applications. One is the need for an appropriate 'persistence
condition', which is needed to ensure that the divisor
filtration and the bundle filtration interact well. The other is the fact
that the divisors involved in the echelon modifications are, in a sense, smaller than the 'obvious' divisors
one could work with. In the situation of curve families with
reducible fibres, this implies working on the \emph{total} space of the family (or
something like it), and constructing modifications
that are not a pullback from the base. This feature
is moreover critical for the universal property of echelon
modifications (Theorem \ref{thm}, (iii)).\par
This note was mostly contained in \cite{grd} originally,
but will
be published separately.
\section{Echelon data and their normal forms}
Our purpose here is to define the notion of \emph{echelon datum}
on a scheme. Roughly speaking, such a datum consists of a vector bundle
$E$ together with a descending chain of full-rank locally free subsheaves
 $E^i$,
such that the degeneracy loci of the inclusion maps $E^i\to E$
form a (usually non-reduced) divisorial chain, more specifically that $E^i$
has minimal generators that have zeros of size $\delta_j, j\leq i$
as sections of $E$, where $\delta_1\leq \delta_2...$ is a suitable
ascending chain of divisors. The definition is as follows.
\begin{defn}\begin{itemize}\item
A \emph{pre-echelon datum} of length $m$, $(E^.,\delta.)$ on a scheme $X$
consists of the following items % (i), (ii), (iii) (resp. (i), (ii)):
\begin{enumerate}
\item a descending chain of locally free sheaves
of the same rank \[E=E^0\supset E^1...\supset E^m;\]
\item a collection of Cartier divisors $\delta_0=0, \delta_1,...,\delta_m$, such that $d\delta_i:=\delta_i-\delta_{i-1}$ is effective;
\end{enumerate}
%\item a descending chain of invertible ideals on $X$
%\eqspl{}{
%\begin{matrix}
%D_1&\to& D_2&...&\to&D_m\\
%\downarrow&&\downarrow&&&\downarrow\\
%\delta_1&\to&\delta_2&...&\to&\delta_m
%\end{matrix}
%\fJ_0=\O\supset\fJ_1=\O(-\delta_1)\supset...\supset\fJ_m=\O(-\delta_m)
%} such that, setting $\Delta\fJ_i:=\fJ_i\fJ_{i-1}\inv$ (as invertible ideal), we have that $E^i\supset\Delta\fJ_i E^{i-1}$ and
which satisfy the following\par {\bf{'persistence condition':}}
\ for all $i$, we have
$E^i\supset E^{i-1}(-d\delta_i)$ and the quotient
 $E^i/E^{i-1}(-d\delta_i)$ maps isomorphically to a
locally free, locally split $\O_{d\delta_i}$-submodule of
$E^{j}\otimes\O_{d\delta_i}$, for all $j<i$. \item An \emph{echelon datum}
$\chi=(E^., \delta_., D.)$ consists of a pre-echelon datum
$(E^.,\delta.)$ plus a collection of Cartier sub-divisors $(D_.\leq
\delta_.)$ such that $dD_i:=D_i-D_{i-1}$ is effective
 and $D_i$ has no components in common with
$D_i^\dag:=\delta_i-D_i$.
\end{itemize}

\begin{comment}
with the following properties:
\begin{enumerate}
\item $E^i\supset\O(-\delta_i)E$;
\item
\ each graded piece
 $G^i:=E^{i-1}/E^i$ is a locally free $\delta_i-\delta_{i-1}$-module
 (where $\delta_0=0$).
\item for all $1\leq i\leq j$, the images of $E^i$ and $E^j$ in
$E|_{\delta_i}$ coincide, i.e.

\[ E^i=E^j+\O(-\delta_i)E\subset E\]

 \end{enumerate}
 \end{comment}
%\end{enumerate}
\end{defn}
\begin{rem}
Note that for all $0\leq j<i$, $E^i$ contains $E^j(\delta_j-\delta_i)$,
so that \[E^i/E^j(\delta_j-\delta_i)\subset E^j\otimes \O_{\delta_i-\delta_j}\] is a well-defined subsheaf,
and we have commutative
\eqsp{
\begin{matrix}
E^i/E^j(\delta_j-\delta_i)&\to&E^i/E^{i-1}(-d\delta_i)\\
\downarrow&&\downarrow\\
E^j\otimes\O_{d\delta_i}&\leftarrow&E^{i-1}\otimes\O_{\delta_i}
\end{matrix}
}

The
persistence condition means that the right arrow is an isomorphism to a locally free and cofree (i.e. split) subsheaf
of its target, and likewise for the composite of the right and bottom arrows. See example \ref{echelon-example-cont} for
motivation for the persistence condition.

\end{rem}
\begin{rem}
An echelon datum is said to be \emph{scalar} if $\delta_i=n_i\delta, D_i=n_iD$ for all $i$ and
fixed effective divisors $\delta, D$. This cases
considered in this paper have this property.\end{rem}
\begin{example}\lbl{echelon-example}
One (scalar)
example to keep in mind is the following: $\pi:X\to B$
is a proper morphism (e.g. a family of curves), $\L$ is a line bundle on $X$, $E=\pi^*(\pi_*(\L))$, $\delta_1=\pi^*(\delta)$,
$D=D_1$ is a component of $\delta_1$, and
\[E^i=\pi^*(\pi_*(\L(-iD))).\qed\]
\end{example}
In analyzing echelon data locally, a useful tool is a normal form
 called \emph{echelon decomposition}, constructed as follows.
Let $t_i$ be an equation for $\delta_i-\delta_{i-1}$. %-\delta_{i-1}$, so that
%$t_1...t_i$ is an equation for $\delta_i$. Also set
%\[ t_{i,j}=\begin{cases} t_i\cdots t_{j+1},~~~ i>j;\\
%1,~~~ i\leq j\end{cases} \]
%this being an equation for $\delta_i-\delta_j$ if $i\geq j$.
%Let $B_0$ be a lift of the image of $E^1\otimes\O_{\delta_1}\to E\otimes\O_{\delta_1}$ to a free submodule of $E$ and let $A_m\subset E$ be a free complement. Then
%\[E^1=B_0\otimes t_1A_m\subset E=B_0\otimes A_m.\]
Using the persistence condition, we can find a free submodule $A_0$ of $E^m$ which maps isomorphically to its
image in $E^{m-1}$, and a free complement to the latter,
\mbox{$A'_1\subset E^{m-1}$,} so that
\[E^m=A_0\oplus t_mA'_1\subset E^{m-1}=A_0\oplus A'_1\]
By assumption, $A_0$ persists in $E^{m-2}$; therefore we can find
free submodules $A_1, A'_2\subset E^{m-2}$ so that
\[E^{m-1}=A_0\oplus A_1\oplus t_{m-2}A'_2\subset E^{m-2}=
A_0\oplus A_1\oplus A'_2,\] therefore as submodule of $E^{m-2}$, we have
\[E^m=A_0\oplus t_mA_1\oplus t_mt_{m-1}A'_{m-2}\]
Continuing in this way and setting finally $A'_m=A_m$ we obtain what
we will call an \emph{ echelon decomposition}
\[ E=A_0\oplus...\oplus A_m\]
such that \eqspl{echelon-decomp}{E^i= A_0\oplus...\oplus
A_{m-i}\oplus t_{i}A_{m-i+1}\oplus ...\oplus t_{i}...t_1 A_m,
i=0,...,m.} In particular,
\[E^m=A_0\oplus t_mA_1\oplus t_mt_{m-1}A_2\oplus...\oplus
t_m...t_1A_m.\]
\begin{example}[Example \ref{echelon-example}, cont]\label{echelon-example-cont}
In the situation of the last example, where $t_i=t=xy$ is an equation for $\delta$ and $y$ is an equation for $D$,
the summand $A_0\oplus...\oplus A_{m-i}$ corresponds to the
image of general sections of $\L(-iD)$ (divisible by
$y^i$ as section of $\L$); $tA_{m-i+1}$ comes from sections
of $\L(-iD)$ divisible by $x$, etc.; $t^iA_m$ comes from
sections of $\L(-iD)$ divisble by $x^i$. Because multiplication by $y$ does not affect linear independence modulo $x$, it is easy to see that
the persistence condition is satisfied.
\end{example}
\section{Echelon modification}
Our purpose is to associate to an echelon datum
$\chi=(E^., \delta., D.)$  a type of
birational modification of the bundle $E=E^0$  which
generalizes the familiar
notion of elementary modification.
This echelon modification will be an \myit{ascending} chain of vector bundles
$$E=E_0\subset E_1...\subset E_m=\Mod(\chi, E),$$ all equal  to $E$ off
 $D_m$. \textdbend: $E_i\neq E^i$. The process works
 iteratively, first producing $E_1$ carrying
  an echelon datum of length $m-1$, etc.
   \par
 %Furthermore, if $E$ is additionally endowed with another echelon datum $\chi'$ compatible with $\chi$, there will be an induced echelon datum on $\Mod(\chi, E)$:
%\[\chi'_m=\Mod(\chi, \chi').\]\par
Let $(E^., \delta., D.)$ denote an echelon datum on $X$. %\overrightarrow{}%, and let
%$D$ be a Cartier divisor on $X$ such that $D_i:=D\cap \delta_i$
%is a Cartier divisor for each $i$ and $D_i-D_{i-1}$ is effective.
%Thus $D_i$ is essentially a sub-divisor, i.e. a sum of components
%of $\delta_i$. There is no loss of generality in assuming
%$D=D_m$.
For $i\geq 1$, set \[dD_i=D_i-D_{i-1},  L_i=\O_{dD_i}(dD_i), i\geq 1, D_0:=0.\]
 $L_i$ is an invertible $\O(dD_i)$-module. Also let\[G^i=\coker(E^i\to E^{i-1}),
 H^i=E^i/\O(-d\delta_i)E^{i-1}\]
which are by assumption  locally free $\O_{\delta_i}$-modules. We
have a locally split exact sequence of  locally free
$\O_{\delta_i}$-modules  \eqspl{}{ \exseq{H^i}{E^{i-1}\otimes
\O_{d\delta_i}}{G^i} }
 This yields a similar exact sequence upon restriction
on $dD_i$. In particular, the kernel of the
natural map $E^{i-1}\to G^i\otimes\O_{dD_i}$ is generated by
$E^i$ and $E(-D_i)$.
For $H^i$, there is another exact sequence
 \eqspl{}{
 \exseq{E^{i-1}(-\delta_i)}{E^{i}}{H^i}
 }
 %Moreover, for all $j\geq i$, the map
 %$E^j\to H^i$ is still surjective.
 In term of a local echelon decomposition \eqref{echelon-decomp},
 we have
 \[H^i=(A_0\oplus...\oplus A_{m-i})\otimes\O_{d\delta_i}\]
\par
To start our ascending chain, define a subsheaf
$E_1\subset E(D_1)$ by the exact diagram
%define a subsheaf  $E_1\subset E(D_1)$ by
%\eqspl{}{
%E_1=E^1(D_1)+E\subset E(D_1).
%}
%Note the exact sequence
\eqspl{}{
\begin{matrix}
0\to & E^1(D_1)&\to& E(D_1)&\to &G^1(D_1)&\to 0\\
     &\downarrow&&||&&\downarrow&\\
0\to & E_1&\to& E(D_1)&\to& G^1\otimes L_1&\to 0.\end{matrix}
}
Because the right vertical map is surjective
with kernel $E\otimes \O_{D_1^\dag},D_1^\dag:=\delta_1-D_1$,
$E_1$ is generated by $E$ and $E^1(D_1)$, i.e.
\eqspl{}{
E_1=E^1(D_1)+E\subset E(D_1).
}
Because $G^1$ is locally free over $\O_{\delta_1}$,
it follows that $G^1\otimes L_1$ is a locally free
$\O_{D_1}$-module, hence
 $E_1$ is an elementary modification of $E$ and in particular is locally
free over $X$. Also, the snake lemma yields exact
\eqspl{E1-I}{\exseq{E^1(D_1)}{E_1}{G^1\otimes\O_{D_1^\dag}}.
}
In particular, $E_1$ is just $E^1(D_1)$ if $D_1=\delta_1$.
We have another the exact diagram
\eqspl{E1-II}{
\begin{matrix}
&0&&0&&&\\
&\downarrow&&\downarrow&&&\\
&E&=&E&&&\\
&\downarrow&&\downarrow&&&\\
0\to&E_1&\to&E(D_1)&\to&G^1\otimes L_1&\to 0\\
&\downarrow&&\downarrow&&\parallel&\\
0\to&H^1\otimes L_1&\to&E\otimes L_1&\to&G^1\otimes L_1&\to 0\\
&\downarrow&&\downarrow&&&\\
&0&&0&&&
\end{matrix}
} In terms of a local decomposition as in \eqref{echelon-decomp}, we can write, where generally $t_i=x_iy_i$ with $y_i$ an equation
for $dD_i$,
\[ E_1= \frac{1}{y_1}(A_0\oplus...\oplus A_{m-1})\oplus A_m.\]
Next, define subsheaves $E^i_1\subset E_1$ by
\[E^i_1=E^{i+1}(\delta_1)\cap E_1\subset E(\delta_1).\]
In particular,
\[E^1_1=E^2(\delta_1)\cap E_1.\]
In terms of a local echelon decomposition, we have
\[E_1^i=\frac{1}{y_1}(A_0\oplus...\oplus A_{m-i-1}
\oplus t_{i+1}A_{m-i}\oplus...\oplus t_{i+1}...t_2A_{m-1})
\oplus t_{i+1}...t_2A_m\]
Then set, analogously to the above,
\[E_2=E_1^1(dD_2)+E_1\subset E_1\otimes dD_2\]
In terms of an echelon decomposition, this is
\[E_2=\frac{1}{y_1y_2}(A_0\oplus...\oplus A_{m-2})
\oplus \frac{1}{y_1}A_{m-1}\oplus A_m\]
Then we have exact
\[\exseq{E_1}{E_2}{H^2\otimes\O(D_2)}\]
Note the natural inclusion
\eqspl{echelon-comparison-eq}{
E^2(D_2)\to E_2
}
In general, we define inductively
\eqspl{}{E_j^i=E_{j-1}^{i+1}(d\delta_j)\cap E_j, i\geq 1\\
E_{j+1}=E_j^1(dD_{j+1})+E_j\subset E_j(dD_{j+1}).}
Again we have an inclusion, $\forall i$,
\[ E^i(D_i)\to E_i.\]

%[[[[[[[[[*********
\begin{comment}
 Next, define $E_2$ as
\eqspl{}{
E_2=E^2(D_2)+E_1\subset E_1(dD_2)\subset E(D_2).
} Generally, if $E_i$ is defined, set
\eqspl{}{
E_{i+1}=E^{i+1}(D_{i+1})+E_i\subset E_i(dD_{i+1})
\subset E(D_{i+1})
}
We have exact sequences
\eqspl{}{
\exseq{E_{i+1}}{E_i(dD_{i+1})}{G^i\otimes L_i},
}
\eqspl{}{
\exseq{E_i}{E_{i+1}}{(E^{i+1}/\O(-dD_{i+1})E^i)\otimes L_{i+1}}.
} With either of these, together with the Riemann-Roch without denominators, divisor version (see \cite{ful}, Example 15.3.4 and \S\ref{asd} below), we can
easily compute, in principle, the Chern classes of all the $E_i.$
Also, it is important to note the injection
\eqspl{E^i-to-E_i}{
E^i(D_i)\to E_i
} whose cokernel is an $\O_{D_i^\dag}$-module; in particular, this injection is a generic isomorphism on $D_i$ and an iso
if $D_i=\delta_i$.
\end{comment}
%****]]]]]]]]]]]]]]]]
\par
In terms of a local echelon decomposition as in \eqref{echelon-decomp}, we can describe the $E_i$ as follows.
\eqspl{}{
E_i=\frac{1}{y_1...y_i}(A_0\oplus...\oplus A_{m-i})\oplus\frac{1}{y_1...y_{i-1}}A_{m-i+1}\oplus...
\oplus A_m,\\
E_m=\frac{1}{y_1...y_m}A_0\oplus\frac{1}{y_1...y_{m-1}}A_1
\oplus...
\oplus\frac{1}{y_1}A_{m-1}\oplus
A_m
}
\begin{rem}
Note that via the various exact sequences (e.g. \eqref{E1-I}, \eqref{E1-II}) above, $K(X)$- group elements,
Chern classes and similar attributes of the
modifications $E_i$ are computable in terms of
similar attributes of the echelon data.\end{rem}
We summarize some of the main properties of elementary modifications as follows. All of them follow directly from the explicit construction and local forms given above. Property (iii), the universal property of echelon modifications, is from our perspective the main raison d'etre for the construction.
\begin{thm}\label{thm}
Let $(E, \chi)$ be an echelon datum on an integral scheme $X$.
\begin{enumerate}\item
$\Mod(\chi, E)$ is a locally free sheaf containing and generically equal to $E$, and
depends functorially on $(\chi, E)$.
\item If $f:Y\to X$ is a dominant morphism from another
integral scheme, then
\[f^*\Mod_X( \chi, E)=\Mod_Y(f^*(\chi), f^*(E)).\]
\item If $\phi:E\to L$ is a map to a line bundle,
such that for each $i$, $\phi(E^i)\subset L(-D_i)$, then
$\phi$ extends to a map
\[\Mod(\chi,\phi):\Mod(\chi,E)\to L.\]
If $E'$ is any locally free sheaf containing and generically equal to $E$ such that $\phi$ extends to
$E'$, then $E'$ is contained in $\Mod(\chi, E)$.

\end{enumerate}
\end{thm}
%[[[[[[[[[
\begin{comment}
There is one notable case where the result is especially simple. Let us say that a subvariety $Z\subset Y$ is \emph{neutral} with respect to an echelon datum $\chi$
as above is $\O_Z(D_i)$ is trivial for all $i$. This can happen if $Z$ is disjoint from $D_i$ or
if $Z\subset D_i$ and the normal bundle $\O_{D_i}(D_i)$ is trivial on $Z$. Because $E$ and $E_1$ admit
 filtrations whose respective associated graded pieces
differ by twisting by some $dD_i$, those gradeds have isomorphic restrictions on $Z$. Therefore
\begin{cor}\lbl{neutral-subvariety-cor}
Notation as above, if $Z$ is neutral with respect to $\chi$
then the total Chern class $c(E_i|_Z)=c(E|_Z)$ for all $i$.
\end{cor}
\par[[[[

 Locally, using an echelon basis $(e_{jk})$, we can describe $E_i$
as follows. Let $y_i$ be a local equation for $dD_i$,
so that $y_i\cdots y_1$ is a local equation for $D_i$.
Then $E_i$ is freely generated by
\eqspl{}{
\frac{1}{y_{\min(i,m-j)}\cdots y_1}e_{jk}, 0\leq j\leq m, 1\leq k\leq k_j
} The coefficient is obtained from
$\frac{t_{i,m-j}}{y_i\cdots y_1}$ by replacing all $x_k$ by 1
(and cancelling $y$'s). Thus $E_i$ is locally free.
\end{comment}
%*******]]]]]]]]
\par
\newsection{Polyechelon data}
In practice, one needs to work with multiple echelon data on a given
bundle. This is feasible provided the data satisfy a reasonable
condition of \myit{transversality}, a strong version of which
follows.
\begin{defn}
Let $(\chi_j=(E^._j,\delta_j., D_j.))$ be a collection of echelon data of respective lengths $m_j$ on
a given bundle $E$. This collection is said to be mutually \emph{transverse}
if for any choice of subset $S\subset \{1,...,m\}$ and
assignment $S\ni j\mapsto i(j)$,\begin{itemize}\item
 the sequence of divisors $(D_{j,i(j)}-D_{j,i(j)-1}: j\in S)$
is regular, i.e. its intersection has codimension $|S|$;
\item for all $i\not\in S$, $(E^k_i\cap \bigcap \limits_{j\in S} E_j^{i(j)}, \delta_{i,k}, D_{i,k}:k=0,...,m_i)$ is an echelon
    datum on $\bigcap\limits_{j\in S} E_j^{i(j)}$.
%\item there is an inclusion
%\[\bigcap\limits_{j\in S} E_j^{i(j)-1}(-\sum\limits_{j\in S}\delta_{j,i(j)})
%\subset \bigcap\limits_{j\in S} E_j^{i(j)}\]
%whose cokernel is $\O_{\sum\limits_{j\in S}\delta_{j,i(j)}}$-locally
%free.
%\item $\forall i',j', J, i(j)$, $q^{i'}_{j'}$ has constant rank on
%$\bigcap\limits_{ j\in J}\ker(q^{i(j)}_j)$.
\end{itemize}
A \emph{polyechelon datum} is a mutually transverse collection of echelon data as above.
\end{defn}
Now a key, albeit elementary, observation
is the following. If $\chi, \chi'$ are transverse echelon data,
with corresponding filtrations $F^.(E) , (F')^.(E)$, then
for each $k$, there is an
 echelon datum $(F^.(E)\cap (F')^k(E), \delta., D.)$
and performing the corresponding echelon modifications leads to a new bundle $\Mod(\chi, (F')^k(E)$. Together these bundles form
 echelon datum
\[\Mod(\chi, \chi')=(\Mod(\chi, (F')^k(E), \delta'_k, D'_k:
k=0,...,m').\]
This is an echelon datum on $\Mod(\chi, E)$.\par
This operation can be iterated: given a transverse collection
of echelon data $(\chi_1,...,\chi_k)$, echelon modifications
 yield a new
transverse collection of echelon data
\[(\chi_{1,2}=
\Mod(\chi_1, \chi_2),...,\chi_{1,k}=\Mod(\chi_1, \chi_k))\]
on $M_1=\Mod (\chi_1, E)$, etc. Iterating, we get an increasing chain

\begin{comment}
Now let $E$ be a bundle with an echelon datum. %%
If $E$ is endowed with another
echelon datum $\chi_1=(E^._1, \delta_{1.}, D_{1.})$ that is
transverse to $\chi$,
then it is easy to check that $\chi$ induces an echelon datum on each
$E^{'i}$ and the  corresponding modification yield an echelon datum on
$\Mod(\chi, E)$ that we denote by $\Mod(\chi, \chi_1)$ . This gives rise to the possibility of defining
poly-echelon modifications, as follows.\par
Let $(\chi_1, ...,\chi_k)$ be a transverse  collection of echelon data on a vector
bundle $E$. Define a sequence of pairs $(M_j, \mu_{j,i}), j\geq i$, where each $(\mu_{j,.})$ is a compatible collection of echelon data
on the vector bundle $M_j$, inductively by:
\begin{itemize}
\item $M_1=\Mod(\chi_1, E), \chi_{1,i}=\Mod(\chi_1, \chi_i), i>1$;
\item $M_{j+1}=\Mod(\chi_{j,j+1}, M_j)$;
\item $\chi_{j+1, i}=\Mod(\chi_{j,j+1}, \chi_{j,i}), i>j+1$.
\end{itemize}
This gives rise again to an increasing sequence of (generically
equal) vector bundles
\end{comment}

\[ E=M_0\subset M_1\subset...\subset M_k\]
and we will call this sequence (or sometimes just its
last member) the \emph{poly-echelon modification} of $E$
with respect to the poly-echelon data $(\chi.)$, denoted
respectively
\[ (M.)=(\Mod.(\chi., E)), M_k=\Mod(\chi., E).\]
The following is elementary
\begin{prop}
Notations as above,\par (i) $(\Mod.(\chi., E))$ is a sequence of locally free sheaves and generic isomorphisms;\par
(ii) $\Mod(\chi.,E)$ is independent of the ordering or the sequence
$(\chi.)$.
\end{prop}
\begin{example} This is the sort of situation we have in mind for echelon data.
Let \[\pi:X\to B\] be a family of nodal curves and $\delta\subset B$ a
boundary component corresponding to a relative separating node $\theta$.
We assume the
boundary family $X_\delta$ splits globally as
\[\lf{X}\cup
\rt{X},\] with the two components having local equations
$x,y$, respectively with $xy=t$ a local equation of $\delta$. Let $L$ be a line bundle on $X$ and
\[(n.)=(n_0=0<n_1<...<n_m)\] be an increasing sequence
of integers with the property that for each $i$,
\[E^i:=\pi_*L(-n_i\rt{X})\] is locally free and its formation
commutes with base-change, as will be the case whenever $R^1\pi_*L(-n_i\rt{X})$ is locally free or equivalently,
$h^1(X_t, L(-n_i\rt{X})\otimes\O_{X_t})$ is independent of $t$.
Let $\rt{D}_i=n_i\rt{X}$.
Then \[\rt{\chi}=(E^., n.\delta, n.\rt{D})\] is
an echelon datum on $X$. Likewise for the analogous $\lf{\chi}$.
Clearly $\rt{\chi}$ and $\lf{\chi}$ are transverse.
We may construct the two echelon modifications
of $(E.)$ \[\rt{M}=\Mod(\rt{\chi}, E), \lf{M}=\Mod(\lf{\chi},E)\] and
\[M_\theta=\Mod(\Mod(\rt{\chi}, \lf{\chi}), \lf{M})=\Mod(\Mod(\lf{\chi}, \rt{\chi}), \rt{M}).\]\par A fundamental point here, which follows from Theorem \ref{thm},
 is that the
natural map $E\to L$ factors through $M_\theta$.
Heuristically that is
because a section $s$ of $E_i$ yields a section of $L$
vanishing to order $n_i$ on $\rt{X}$, hence
 $s/y^{n_i}$ still yields a
regular section of $L$.
\end{example}
\vfill\eject
\nocite{hart}\nocite{okonek}\nocite{eisenbud-harris}
\bibliographystyle{amsplain}
\bibliography{mybib}
\end{document}